\documentclass[a4paper,11pt,reqno]{amsart}

\usepackage{latexsym}
\usepackage{amsthm,amssymb,amsmath,amsopn}  
\usepackage{graphics}
\usepackage[all]{xypic}

\newtheorem{theorem}{Theorem}[section]
\newtheorem{lemma}[theorem]{Lemma}

\newcommand{\sfrac}[2]{{\scriptstyle{\frac{#1}{#2}}}}

\newcommand{\e}{\mathrm{e}}

\newcommand{\N}{\mathbf{N}}

\newcommand{\R}{\mathbf{R}}

\newcommand{\wb}{$\circ$}
\newcommand{\bb}{$\bullet$}

\renewcommand{\epsilon}{\varepsilon}
\renewcommand{\theta}{\vartheta}

\begin{document}

\begin{abstract}
In 2003,
Mar{\'o}ti showed that one could use the machinery of 
$\ell$-cores
and $\ell$-quotients of partitions to establish lower 
bounds for~$p(n)$,
the number of partitions of~$n$.
In this paper we explore these ideas in the case $\ell=2$,
using them to give a largely combinatorial proof
of an effective upper bound on $p(n)$, and to prove asymptotic
formulae for the number of self-conjugate partitions,
and the number of partitions with distinct parts.
In a further application we give a combinatorial
proof of an identity originally due to Gauss.
\end{abstract}

\title{Counting partitions on the abacus}
\author{Mark Wildon}
\date{August 3, 2006. }
\email{M.J.Wildon@swansea.ac.uk.}
\maketitle

\section{Introduction}
A \emph{partition} of a number $n \in \mathbf{N}_0$ is a
sequence $(\lambda_1, \ldots, \lambda_k)$ such that 
$\lambda_1 \ge \lambda_2 \ge \ldots \ge \lambda_k \ge 1$
and $\lambda_1 + \ldots + \lambda_k = n$. 
To indicate that $\lambda$ is a partition of $n$ we write $|\lambda| = n$.
Let $p(n)$ be the number of partitions of~$n$. 
By our definition, $\emptyset$ is the unique partition of $0$.

The problem of finding asymptotic estimates for $p(n)$
was introduced, and one might reasonably say, solved, by the
famous paper of Hardy--Ramanujan~\cite{HardyRamanujan}.
As well as their  exact formula for~$p(n)$, proved as the 
original application of the circle-method,
their paper also contains several weaker results proved using
less sophisticated analytic methods; for example, that
\[ \log p(n) \sim c\sqrt{n} \quad\text{as $n \rightarrow \infty$} \]
where $c = 2\sqrt{\pi^{2}/6}$. The only time they refer to using
combinatorial arguments to obtain lower and upper bounds for~$p(n)$ 
is in a footnote in~\S 2.1, which suggests the strongest such bounds 
they could prove were of the form
\[ 2^{A\sqrt{n}} < p(n) < n^{B\sqrt{n}}. \]

Our first aim is this paper is to improve on this. In \S 2
we introduce  G.~D.~James' abacus notation for partitions and use it
to establish an important bijection on partitions (see \S 2.2). 
In \S 3 we use this bijection to show that for every $\epsilon > 0$ 
there exists a constant $A(\epsilon)$ such that
\begin{equation}\label{upper} 
\log p(n) \le A(\epsilon) n^{1/2 + \epsilon} 
\quad \text{for all $n \ge 0$}.
\end{equation}

Despite the appearance of an `$\epsilon$' in this inequality, the 
proof of \eqref{upper} is almost entirely combinatorial. The closest competing
result in the literature is due to Erd{\H{o}}s, who gives at the start 
of \cite{ErdosPartitions} a short elementary proof that
\begin{equation}
\label{Erdosbound} 
	\log p(n) \le c n^{1/2} \quad\text{for all $n \ge 0$},
\end{equation}
where, as before, $c = 2\sqrt{\pi^{2}/6}$. 
Erd{\H{o}}s' result is of course significantly
sharper than~\eqref{upper}, but this improvement comes at the cost of 
a greater contribution from analysis. For instance, the
reader may already have correctly guessed that  the identity
$\sum_{n=1}^\infty 1/n^{2} = \pi^{2}/6$ is needed.

In this context, we mention that Peter Neumann (Oxford Algebra 
Kinderseminar, January 2003) has asked whether there is a combinatorial
proof that there exists a constant $C$ such that 
$\log p(n) \le C\sqrt{n}$ for all $n \ge 0$. While \eqref{upper} 
comes tantalisingly close, the ideas behind it do not appear quite 
powerful enough to solve his problem.

Our second aim is to show that the partition bijection given in \S 2.2
may be used to solve several other problems in the enumeration
of partitions. In our first application we give
a very short proof of an identity originally due
to Gauss (see~\cite[end of \S 3]{Gauss}):
\begin{equation}
\label{Gauss} 
\frac{(1-x^2)(1-x^4)(1-x^6) \ldots}{(1-x)(1-x^3)(1-x^5) \ldots} = 
\sum_{m \ge 0}
x^{\frac{m(m+1)}{2}}.
\end{equation}

Our main application  concerns a result of Erd{\H{o}}s.
In the main part of~\cite{ErdosPartitions}, Erd{\H{o}}s uses the same
underlying idea as he used to prove \eqref{Erdosbound}, and some 
more lengthy (but still elementary) analysis,
to prove that there exists a constant~$b \in \mathbf{R}$ such that
\begin{equation}
\label{asympt} p(n) \sim \frac{\e^{c \sqrt{n}}}{bn} \quad\text{as $n 
\rightarrow \infty$}.
\end{equation}
It was already known from the work of Hardy and Ramanujan 
\cite{HardyRamanujan} that \hbox{$b = 4\sqrt{3}$}, but
Erd{\H{o}}s' methods were not sufficiently powerful to prove this. 
Later however, Newman  \cite{Newman} gave an elementary proof that if 
\eqref{asympt} holds for \emph{some} constant $b$, then necessarily 
$b =4\sqrt{3}$. In \S 4 we use our bijection to
give an alternative proof of this result. In length
it is comparable with Newman's, but it requires considerably
less knowledge of analysis, and the motivation is perhaps more apparent.

In \S 5 we conclude by showing how the abacus
may be applied to give exact asymptotic
estimates for two other classes of partitions, namely 
self-conjugate partitions, and partitions with distinct parts. 

\section{A partition bijection}
Let $\lambda$ be a partition of $n$. A $2$-\emph{hook} in $\lambda$ 
consists of two adjacent nodes in the Young diagram of $\lambda$ whose 
removal leaves the diagram of a partition of $n-2$. By repeatedly 
removing  $2$-hooks from $\lambda$ one obtains the \emph{$2$-core} 
of~$\lambda$. For example, the $2$-core of $(6,3,3,1)$ is $(2,1)$, 
obtained after removing the five~$2$-hooks indicated in the diagram below.

\begin{figure}[b]
\scalebox{0.5}[0.5]{\includegraphics{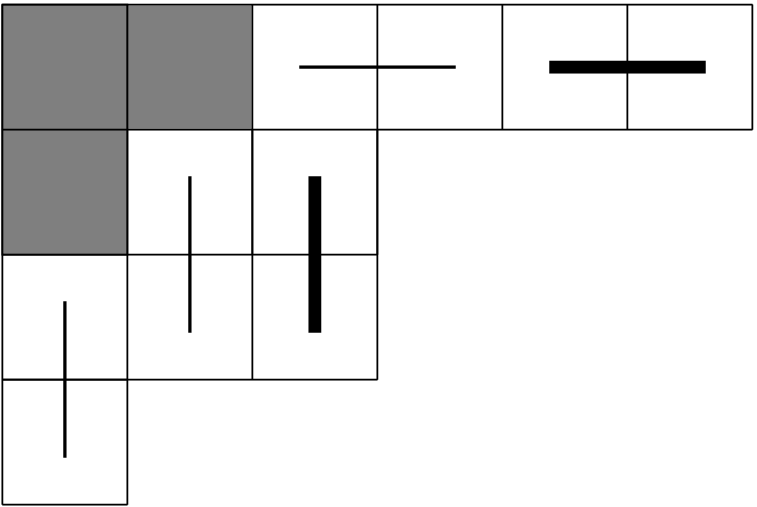}}
\caption{The $2$-core of $(6,3,3,1)$.}
\end{figure} 

In this diagram, two of the three 
$2$-hooks of the original partition are shown by bold lines.
The remaining lines are $2$-hooks of partitions obtained
\emph{en route} to the $2$-core.

It is not obvious that the $2$-core of a partition is independent
of the manner in which one removes its $2$-hooks. This, and much
else besides, can be seen very clearly if we represent partitions 
using G.~D.~James' abacus (see \cite[Ch.~2, 26--27]{JK}). 
To make this article self-contained we briefly recall how
to operate this piece of apparatus. (We describe
only  the `binary' $2$-runner abacus.)

\subsection{}
Let $\lambda$ be a partition of $n$. Starting in the southwest corner
of the Young diagram of $\lambda$ walk along its rim, heading 
towards the northeast corner. For each  step right, put a space, 
indicated~\wb, and for each step up,
put a  bead, indicated~\bb. For example, the partition
$(6,3,3,1)$ shown above has the sequence 
\wb\bb\wb\wb\bb\bb\wb\wb\wb\bb.
It will be seen shortly that it is useful 
to allow such a sequence to begin with any number of 
beads, and to finish with any number of spaces --- these must be 
stripped off before the partition is recovered. 

One then arranges the sequence in two columns, known as the \emph{runners}
of the abacus. For instance $(6,3,3,1)$ is represented by
\[ \begin{matrix} $\wb$ & $\bb$ \\ $\wb$ & $\wb$ \\ $\bb$ & $\bb$ \\ 
$\wb$ & $\wb$ \\ 
$\wb$ & $\bb$ 
\end{matrix}\ . \]

We call this an \emph{abacus} display for the partition $\lambda$. 
One sees that $2$-hooks in~$\lambda$ correspond to  beads with a 
space above them, and that an abacus display for the partition 
obtained by removing a given $2$-hook is obtained by sliding
the corresponding bead  one space up its runner.
(It is at this point that our convention about
initial beads and final spaces is needed.) The $2$-core of $\lambda$
is obtained by pushing all the beads up as far as they will go.

To reconstruct the partition $\lambda$ from an abacus display for its $2$-core
we need to know the number of positions down each bead must be moved. 
For 
example, given the abacus display 
\[ \begin{matrix} $\bb$ & $\bb$ \\ $\wb$ & $\bb$ \\ $\wb$ & $\bb$
\\ $\wb$ & $\wb$ \\ $\wb$ & $\wb$
\end{matrix}
\]
for the $2$-core of $(6,3,3,1)$, we recover the original partition 
by moving the bead on the first runner down by~$2$ positions, 
and on the second runner moving the lowest bead down by~$2$ positions 
and the next lowest down by~$1$ position. 
We record this information in a pair of partitions, $(\mu,\nu)$,
known as the $2$-\emph{quotient} of $\lambda$; here $(6,3,1,1)$ has 
$2$-quotient $((2),(2,1))$. 
Note that the total number of nodes in $\mu$ and $\nu$
is equal to the total number of $2$-hooks removed, so
$|\lambda| = |\gamma| + 2(|\mu| + |\nu|)$ where $\gamma$
is the $2$-core of $\lambda$.

\subsection{}
Given a partition $\lambda$ we may, by adding
a bead at the start if necessary, choose an abacus display
for $\lambda$ so that the $2$-core of $\lambda$ has at least
as many beads on the second runner as the first. With
this convention to fix the order of the partitions in the $2$-quotient
of $\lambda$, the correspondence just described between partitions 
and pairs $(\gamma,(\mu,\nu))$ of $2$-cores and $2$-quotients
is bijective. Moreover, it is easy to see that any $2$-core
is of the form $(m,m-1,\ldots,1)$ for some $m \in \N_0$. 
Therefore we have shown that, if $t(n)$ is the number
of pairs $(\mu,\nu)$ where $\mu$ and $\nu$ are partitions
with $|\mu|+|\nu| = n$, then
\begin{equation}
\label{bij}
p(n) = \sum_r t\left( \frac{n-r(r+1)/2}{2}\right)
\end{equation}
where the sum is over all the non-negative integers
$r$ such that $n-r(r+1)/2$ is a non-negative even integer. For example,
\[ p(10) = t(5) + t(2) + t(0). \]
As it is easy to calculate $t(n)$ given $p(m)$ for $m \le n$
using the formula $t(n) = \sum_{m=0}^n p(m)p(n-m)$,
\eqref{bij} gives us a recurrence relation for the values of
the partition function. We exploit this recurrence in \S 3 and \S 4 below.

\subsection{}We now prove Gauss' identity.
Let $P(x)$ be the ordinary generating function for $p$,  
\begin{equation}\label{Pgenfunc}
P(x) = p(0) + p(1)x + p(2)x^2 + \ldots 
= \prod_{n=1}^\infty \frac{1}{1-x^{n}} .
\end{equation}
Since each node in the $2$-quotient of a partition
contributes two nodes to the original partition, our
bijection shows that
\[ P(x) = P(x^2)^2 \sum_{m \ge 0} x^{\frac{m(m+1)}{2}} . \]
By \eqref{Pgenfunc}, this is equivalent to Gauss' identity \eqref{Gauss}.

It should be said that this is not the first combinatorial
proof of Gauss' identity. In \cite{Andrews}, Andrews gives a
bijective proof (in the style of Franklin's proof
of Euler's Pentagonal Number Theorem) that: 
\begin{equation*}
\frac{(1-q^2)(1-q^4)(1-q^6) \ldots}{(1+q)(1+q^3)(1+q^5) 
\dots} = 1 + \sum_{n=1}^\infty (-1)^{n} q^{n(2n-1)}(1+q^{2n}).
\end{equation*}
Replacing $q$ with $-x$ gives \eqref{Gauss}.
Another strategy is to first give a combinatorial
proof of Jacobi's triple product identity, as done for 
example in \hbox{\cite[\S 6.2]{Pak}}; from Jacobi's identity, 
it only takes a simple specialisation
to get~\eqref{Gauss}.                  

\subsection{}
We conclude this section by remarking that the bijection we have 
described may be thought as being between the set of all partitions,
and the set of all oriented binary trees with vertices labelled
by members of $\N_0$, with the property that if two leaves
are connected to the same node, they are not both labelled by 0.
For example, the $2$-core and $2$-quotient of $(6,3,3,1)$
may be represented by the diagram below.

\hfil
\scalebox{0.85}{
\xymatrix{
      & 2\ar@{-}[dl]\ar@{-}[dr] & \\
  (2) &                         & (2,1)
}}\hfill

\smallskip
\noindent
Here the upper `$2$' records the number of rows (or columns)
in the $2$-core of $(6,3,3,1)$.
Iterating the process, we find that $(6,3,3,1)$
is represented by the oriented binary tree
shown below.

\hfil
\scalebox{0.85}{
\xymatrix{ 
  &  &  2\ar@{-}[dl]\ar@{-}[dr] \\
  &  0\ar@{-}[dl]\ar@{-}[dr] &  & 2 & \\
0 &  &  1 \\
}}\hfill

\medskip

The corresponding formulation of Gauss' identity is the attractive
\[ \prod_{n=1}^\infty \frac{1}{1-x^n} = 
\prod_{m=0}^\infty \left( \sum_{r=0}^\infty x^{\frac{r(r+1)}{2}2^{m}} 
\right)^{\!2^{m}}.\]

\section{Lower and upper bounds for $p(n)$}
By counting only partitions of $4m$ with $2$-core $\emptyset$ and $2$-quotient
$(\mu, \nu)$ where both $\mu$ and $\nu$ are partitions of $m$, we
obtain the inequality
$p(4m) \ge p(m)^2$. 
An immediate consequence is that
\[ p(4^r) \ge p(4)^{2^{r-1}} = 5^{2^{{r-1}}} \quad\text{for all $r \ge 1$}. \]
Given  $n \in \mathbf{N}$, we may choose a power of $4$, say $4^r$, such
that $n/4 < 4^r \le n$. Since $p$ is an
increasing function, $p(n) \ge p(4^r)$, and so
\[ \log p(n) \ge \log p(4^r) \ge \frac{\log 5}{2}  2^{r} \ge \frac{\log 5}{4} 
\sqrt{n} \quad\text{for all $n \ge 4$}. \] 
It seems striking that such a simple argument gives a lower bound which,
while far from optimal, is still of the right asymptotic form.

Mar{\'o}ti (see \cite[Corollary 3.1]{MarotiLower})  has shown
that, with a little more analysis, and the use of a computer to check
small cases, one can obtain the significantly stronger lower bound 
\[ p(n) \ge \frac{\e^{2 \sqrt{n}}}{14} \quad\text{for all $n \ge 1$}. \]

The next lemma shows how our bijection may be applied to give an upper bound.

\begin{lemma}\label{lemma:fbound}
Suppose that $f : \R_{\ge 0} \rightarrow \R_{\ge 0}$ is an increasing function 
such that $f(0) \ge 1$ and 
\[ \sum_{0 \le s \le n/2} f\left(\frac{n}{2} - s\right) 
f\left(s \right)
\le \frac{f(n)}{n} \quad\text{for all $n\in\mathbf{N}$.}\]
 Then $p(n) \le f(n)$ for all $n \ge 1$.
\end{lemma}

\begin{proof}
We work by induction on $n$. Using \eqref{bij} we have
\[ p(n) \le n t(n/2) = n \sum_{0 \le s \le n/2} 
p\left(\frac{n}{2} - s\right)p(s) \]
where we adopt the convention that $p(x) = p(\left\lfloor x
\right\rfloor)$, and similarly for $t(x)$. Hence
\[ p(n) \le n \sum_{0 \le s \le n/2} 
f\left(\frac{n}{2} - s\right) f(s) \le f(n) \]
which completes the inductive step.
\end{proof}

The reader will see that we have been very crude in taking $n$
as an upper bound for the number of summands in \eqref{bij} ---
there appears to be no advantage in being any more precise at 
this point.

\begin{theorem}
Given any $\epsilon > 0$ there is a constant $A(\epsilon)$ depending
only on~$\epsilon$ such that
\[ p(n) \le \e^{A(\epsilon) n\raisebox{3pt}{\tiny $\sfrac{1}{2} \!+ \!\epsilon$}}\quad
\text{for all $n \ge 0$}.\]
\end{theorem}

\begin{proof}
Without loss of generality, assume that $\epsilon < \sfrac{1}{2}$.
Let $f(x) = \e^{A x^{\sfrac{1}{2} + \epsilon}}$.
We shall prove that when $A$ is large enough,
$f$ satisfies the conditions of Lemma~\ref{lemma:fbound}.
Setting $\beta = \frac{1}{2} + \epsilon$, we have
\[ \sum_{0 \le s \le n/2} f\left(n/2 - s\right) 
f\left(s\right)
= \sum_{0 \le s \le n/2} \e^{A\left( n/2 - s \right)^\beta
+ A  s^\beta} \le  \frac{n+1}{2} \,\e^{A n^\beta 4^{-\epsilon}} \] 
where we used the inequality
\[ x^\beta + y^\beta \le 2 \left(\frac{x+y}{2}\right)^\beta \]
which is valid for $0 \le \beta \le 1$ and $x,y \ge 0$.
Therefore the only condition we must satisfy is that
\[ \e^{A  n^\beta 4^{-\epsilon}} \le \frac{2}{n(n+1)}\, \e^{A 
n^\beta} \quad\text{for all $n \ge 1$}, \]
or, equivalently, that
\[ \e^{A (1-4^{-\epsilon}) n^\beta} \ge \frac{n(n+1)}{2}
\quad\text{for all $n  \ge 1$}.\]
Clearly this will hold provided $A$ is sufficient large. 
\end{proof}

\section{A new evaluation of the constant in Erd{\H{o}}s' formula}

In this section we assume Erd{\H{o}}s' result \eqref{asympt} that
$p(n) \sim \e^{c\sqrt{n}}/bn$ for some constant $b \in \mathbf{R}$, and prove
that $b = 4\sqrt{3}$. We shall also need Erd{\H{o}}s' 
upper bound~\eqref{Erdosbound} that $p(n) \le \e^{c\sqrt{n}}$
for all $n \ge 1$. (Here, as always, $c = 2\sqrt{\pi^2/6}$.)

Recall that $t(n)$ denotes the number of pairs $(\mu,\nu)$
of partitions such that $|\mu| + |\nu| = n$. We first show that it follows
from Erd{\H{o}}s' result that 
\begin{equation}
\label{pairsasympt} 
t(n) \sim \frac{\e^{c\sqrt{2n}}}{n^{5/4}} \frac{2^23^{1/4}}{b^{2}} \quad
\text{as $n \rightarrow \infty$.}
\end{equation}
We then use the bijection described in \S 2.2 to deduce that
\begin{equation}
\label{lim}
p(32m) \sim \frac{\e^{c\sqrt{32m}}}{32m} \frac{2^23^{1/2}}{b^{2}} \quad
\text{as $n \rightarrow \infty$.}
\end{equation}
Comparing this with $\eqref{asympt}$ shows that $b = 4\sqrt{3}$.
(The only reason for taking~$32m$ rather than a general~$n$ in \eqref{lim} 
is that this choice lead to some simplifications in 
the expressions we encounter;
there is of course no need to prove anything stronger.)

\subsection{}
All the results in this section
and \S 5 are proved using the following  simple lemma,
whose proof requires only the basic integral
\[ \int_0^\infty \e^{-\gamma x^2} \:\mathrm{d}x= \sqrt{\frac{\pi}{4\gamma}} \quad
\text{if $\gamma > 0$.}\]

\begin{lemma}
\label{lemma:gen}
If $\alpha, \beta, \gamma, \theta > 0$ then
\[ \sum_{r=0}^{\alpha m^{\beta + \theta}} \e^{-\gamma r^{2}/m^{2\beta}} \sim
\sqrt{\frac{\pi}{4\gamma}} \,m^\beta \quad\text{as $m\rightarrow \infty$}.\]
\end{lemma}

\begin{proof}
Let $S_m$ stand for the sum on the left-hand-side. Since
\[ S_m -1 \le \int_0^{\alpha m^{\beta+\theta}} \e^{-\gamma x^2/m^{2\beta}} \mathrm{d}x
\le S_m, \]
it is sufficient to prove that
\[ J_m = \frac{1}{m^\beta} \int_0^{\alpha m^{\beta+\theta}} \e^{-\gamma x^2/m^{2\beta}} 
\mathrm{d}x
\rightarrow \sqrt{\frac{\pi}{4\gamma}} \quad\text{as $m\rightarrow\infty$}.\]
This is not hard. On the one hand,
\[ J_m \le \frac{1}{m^\beta}\int_0^\infty \e^{-\gamma x^2/m^{2\beta}} \mathrm{d}x = 
\sqrt{\frac{\pi}{4\gamma}}
\]
and on the other,
\[
\sqrt{\frac{\pi}{4\gamma}} - J_m =
\frac{1}{m^\beta} \!
\int_{\alpha m^{\beta+\theta}}^\infty \!\e^{-\gamma x^2/m^{2\beta}} 
\mathrm{d}x =
\frac{1}{m^\beta} \! \int_{\alpha m^\theta}^\infty \! \e^{-\gamma y^2} 
 m^\beta \mathrm{d}y  
=\int_{\alpha m^\theta}^\infty \!\e^{-\gamma y^2} \mathrm{d}y
\]
which tends to $0$ as $m \rightarrow \infty$. Hence $J_m$ has 
the claimed limit. 
\end{proof}

\subsection{} We start by proving \eqref{pairsasympt}.
If \eqref{pairsasympt} holds when $n$ is even then, by a simple 
sandwiching argument, it must hold for all $n$. We may therefore 
assume that $n = 2m$, in which case \eqref{pairsasympt} becomes
\begin{equation}
\label{qeq} 
t(2m) \sim \frac{\e^{2c\sqrt{m}}}{m^{5/4}} \frac{2^{3/4}3^{1/4}}{b^{2}}
\quad\text{as $n \rightarrow \infty$.} 
\end{equation}
Let $\epsilon > 0$ be given. Choose $N$ so that
\[ 1-\epsilon < \frac{p(n)}{e^{c\sqrt{n}}/bn} < 1 + \epsilon \quad
\text{for all $n \ge N$}. \]

We first obtain an upper estimate for $t(2m)$.
For $0 < \alpha < 1$ set
\begin{align*}
M_{\alpha}(m) &= \sum_{0 \le r \le \alpha m}
p(m+r)p(m-r), \\ 
R_{\alpha}(m) &= \sum_{\alpha m < r \le m}p(m+r)p(m-r).
\end{align*}
so we have $t(2m)/2 \le M_{\alpha}(m) + R_{\alpha}(m)$. 
We shall see that as $m$
tends to infinity, the contribution of $M_{\alpha}(m)$ to~$t(2m)$ dominates, 
no matter what the choice of~$\alpha$. In fact, it follows
from Erd{\H{o}}s' bound~\eqref{Erdosbound} and the inequality
\[ \sqrt{1+x} + \sqrt{1-x} \le 2 - \frac{x^2}{4} \quad\text{if $0 \le x \le 1$} \]
that
\[
R_{\alpha}(m) \le  \sum_{ \alpha m < r \le m }\e^{c\sqrt{m}\left(\sqrt{1+r/m}
+\sqrt{1-r/m}\right)} \le (1-\alpha)m e^{2c\sqrt{m} - 
cr^2/4m^{3/2}},
\] 
and so
\[ \frac{R_{\alpha}(m)}{\e^{2c\sqrt{m}}/m^{5/4}} \le 
(1-\alpha)m^{9/4}
\e^{-c\alpha^{2}\sqrt{m}/4} \rightarrow 0 \quad\text{as $m\rightarrow \infty$}.\]
We now look at the main contribution. Provided $\alpha < 1/2$ and
$m \ge 2N$, the
smallest $n$ for which $p(n)$ appears in $M_{\alpha}(m)$ is at least
$N$. Hence if
these conditions hold,
\begin{align*} 
\frac{M_{\alpha}(m)}{(1+\epsilon)^{2}} &\le
\sum_{r = 0}^{\alpha m}
\frac{\e^{c\sqrt{m}\left(\sqrt{1+r/m}+\sqrt{1-r/m}\right)}}
{b^2(m^2-r^2)} \\
& \le \frac{\e^{2c\sqrt{m}}}
{b^2m^2(1-\alpha^2)}
\sum_{r = 0}^{\alpha m} \e^{c\sqrt{m}
\left(\sqrt{1+r/m}+\sqrt{1-r/m} - 2\right)} \\
& \le \frac{\e^{2c\sqrt{m}}}{b^2m^2(1-\alpha^2)}
\sum_{r = 0}^{\alpha m} \e^{-c r^{2}/4m^{3/2}}.
\end{align*}
It now
follows from Lemma \ref{lemma:gen}, applied with $\beta = 3/4$, $\theta = 1/4$, that
\[ \limsup \frac{M_\alpha(m)}{\e^{2c\sqrt{m}}/m^{5/4}} \le
\frac{(1+\epsilon)^2}{b^2(1-\alpha^2)} \sqrt{\frac{\pi}{c}} =
\frac{(1+\epsilon)^2 3^{1/4}}{2^{1/4}b^2(1-\alpha^2)}.\]
Finally, let $\alpha, \epsilon \rightarrow 0$ to get
\[ \limsup \frac{t(2m)}{\e^{2c\sqrt{m}}/m^{5/4}} \le
\frac{ 2^{3/4}3^{1/4}}{b^2}. \]

We now obtain a lower estimate for $t(2m)$.
By taking $\alpha$ sufficiently small we may ensure that
\[ \sqrt{1+x} + \sqrt{1-x} \ge 2 - \frac{(1+\epsilon)x^2}{4} \quad
\text{if $0 \le x < \alpha$.} \]
For this choice of $\alpha$ we have, by the same manipulations as before,
\[ \frac{t(2m)}{e^{2c\sqrt{m}}/m^{5/4}} 
+ \frac{p(m)^2}{e^{2c\sqrt{m}}/m^{5/4}}
\ge \frac{2}{b^2m^{3/4}} 
\sum_{r=0}^{\alpha m} \e^{-cr^{2}(1+\epsilon)/4m^{3/2}}
. \]
The second term on the left-hand-side compensates for the double 
counting of the contribution of $p(m)p(m)$ to $t(2m)$; as
$p(m) \le (1+\epsilon)\e^{c\sqrt{m}}/bm$ for $m \ge N$, this term vanishes in the limit. As before, it follows from Lemma~\ref{lemma:gen} that
\[ \liminf \frac{t(2m)}{e^{2c\sqrt{m}}/m^{5/4}} \ge 
\frac{1}{b^2} \sqrt{\frac{\pi}{c(1+\epsilon)}} = 
\frac{2^{3/4}3^{1/4}}{b^2\sqrt{1+\epsilon}}. \]
Now let $\epsilon \rightarrow 0$. This completes the proof
of~\eqref{qeq}.

\subsection{}
It follows from \eqref{bij} that 
\[ p(32m) = \sum_r t\left(16m-\frac{r(r+1)}{4} \right) \]
where the sum is over all non-negative integers $r$ 
such that $16m-r(r+1)/4$ is a non-negative even number. 
We may estimate this sum using the same ideas as \S 4.2.

We start with a lower estimate.
The contribution from those $r \equiv 0$ mod~$4$ is at least
the contribution from those $r \equiv 3$ mod~$4$, so we have
\[ p(32m) \ge 2\sum_{s=1}^{\alpha \sqrt{m}} t(16m-(4s-1)s)\]
for any $\alpha$ such that $0 < \alpha < 1$. 

Let $\epsilon > 0$ be given. By \eqref{pairsasympt} we may choose $N$ so that
\[ t(n) > (1 - \epsilon){\frac{\e^{c\sqrt{2n}}}{n^{5/4}} 
\frac{2^23^{1/4}}{b^{2}}}  \quad \text{for all $n \ge N$}. \]
Take $m \ge N$. Choose $\alpha$ sufficiently small that
\[ \sqrt{1-x} \ge 1 - \frac{(1+\epsilon)x}{2} 
\quad\text{if $0 \le x \le \frac{\alpha^2}{4}$.} \] 
For this choice of $\alpha$ we have
\begin{align*}
\frac{p(32m)}{1-\epsilon} &\ge 2\sum_{s=1}^{\alpha \sqrt{m}}
\frac{2^23^{1/4}}{b^2}
\frac{\e^{c\sqrt{32m-2s(4s-1)}}}{(16m-s(4s-1))^{5/4}} \\
&\ge  \frac{2^33^{1/4}}{b^2}  \sum_{s=1}^{\alpha \sqrt{m}}
\frac{\e^{c\sqrt{32m} \sqrt{1-\frac{s^2}{4m}}}}{(16m-s(4s-1))^{5/4}} \\
&\ge \frac{\e^{c\sqrt{32m}}}{32m}
\frac{2^33^{1/4}}{b^2m^{5/4}}
\sum_{s=1}^{\alpha \sqrt{m}} \e^{-(1+\epsilon)cs^2/\sqrt{2m}}.
\end{align*}
It now follows from Lemma \ref{lemma:gen}, applied with
$\beta = 1/4$, $\theta = 1/4$, that
\[ \liminf \frac{p(32m)}{\e^{c\sqrt{32m}}/32m}
\ge \frac{(1-\epsilon)}{\sqrt{1+\epsilon}}
\frac{2^33^{1/4}}{b^2}
\sqrt{\frac{\pi\sqrt{2}}{4c}} 
= \frac{(1-\epsilon)}{\sqrt{1+\epsilon}} \frac{2^{2} 3^{1/2}}{b^2}. \]
Now let $\epsilon \rightarrow 0$. 

To complete the proof of \eqref{lim} we must also show that
\[ \limsup  \frac{p(32m)}{\e^{c\sqrt{32m}}/32m} \le \frac{2^23^{1/2}}{b^2}.\]
In order to apply the method used in \S4.2 we need a uniform upper
bound for $t(n)$; the relatively crude estimate
\[ t(n) \le n p(n/2)^2 \le n\e^{c\sqrt{2n}} \quad\text{for $n\ge 1$}\]
given by \eqref{Erdosbound}
is sufficient for this purpose. 
As no other new ideas are needed
we shall omit the remaining details of the proof.

\section{Asymptotic formulae for two special types of partition}
In this final section we give asymptotic formulae
for the number of self-conjugate partitions, and for the
number of partitions with distinct parts. 
As corollaries, we get precise estimates
for the proportions of such partitions. We omit the details of the analytic
arguments required, as they are of the same nature as those
already seen in \S 4.

\subsection{Self-conjugate partitions}

Recall that if $\lambda = (\lambda_1, \ldots, \lambda_k)$ is a partition
of $n$ then the conjugate partition $\lambda'$ is defined by
$\lambda'_j = |\{ i : \lambda_i \ge j \}|$; the diagram
of $\lambda'$ is obtained from that of $\lambda$ by reflecting
it in its leading diagonal. Let $s(n)$ be the
number of self-conjugate partitions of $n$.

Given an abacus display
for a partition $\lambda$, one obtains a display for the conjugate partition
$\lambda'$ by reading the sequence of beads and spaces from 
right to left,  and then replacing \bb~with~\wb\ and \emph{vice versa}.
For example, $(6,3,3,1)$ and $(6,3,3,1)' = (4,3,3,1,1,1)$ have the
abacus displays:
\[ \begin{matrix} $\wb$ & $\bb$ \\ $\wb$ & $\wb$ \\ $\bb$ & $\bb$ \\ 
$\wb$ & $\wb$ \\ 
$\wb$ & $\bb$ 
\end{matrix}\ ,\qquad\qquad
\begin{matrix} $\wb$ & $\bb$ \\ $\bb$ & $\bb$ \\ $\wb$ & $\wb$ \\ 
$\bb$ & $\bb$ \\ 
$\wb$ & $\bb$ 
\end{matrix} \ 
\]
respectively.

If $\lambda$ is self-conjugate then the sequence of spaces and
beads corresponding to $\lambda$ has an even number of entries. 
It follows that if $\lambda$ has $2$-quotient $(\mu,\nu)$ then
$\nu = \mu^\prime$. This  gives us a bijection
between self-conjugate partitions of~$n$ and
pairs $(\gamma,\mu)$ where $\gamma$ is a $2$-core (and so of a known shape)
and $\mu$ is a partition of $(n-|\gamma|)/4$. (This bijection is well known,
see for instance \cite[Theorem 6.2.2]{JK}. We leave it to the
reader to work out the connection with the labelled binary trees of \S 2.4.)

The corresponding numerical result is
\[ s(n) = \sum p\left(\frac{n-r(r+1)/2}{4}\right) \]
where the sum is over all the non-negative integers $r$ such that $n-r(r+1)/2$ is 
a non-negative integer divisible by $4$. By the same techniques
used in \S 4.3, one finds that
\[ s(n) \sim  \frac{\e^{\frac{c}{2} \sqrt{n}}}{{2^{7/4}3^{1/4}n^{3/4}}}
\quad\text{as $n\rightarrow \infty$}. \]
It follows that the proportion of self-conjugate partitions
satisfies
\begin{equation} 
\frac{s(n)}{p(n)} \sim (6n)^{1/4} \e^{-\frac{c\sqrt{n}}{2}}\quad
\text{as $n\rightarrow \infty$}.
\end{equation}
Thus of the $p(n)$ partitions of $n$, 
approximately $\sqrt{p(n)}$ are self-conjugate.

\subsection{Partitions with distinct parts}
Let $q(n)$ be the number of partitions of $n$ with distinct parts.
Here the numerical result we require is
\[ q(n) = \sum p\left(\frac{n-r(r+1)/2}{2}\right) \]
where the sum is over all the non-negative integers $r$ such that $n-r(r+1)/2$ is a non-negative even integer. 
For a combinatorial proof of this result, see Proposition~5.2 in~\cite{KOR}. 
As an alternative, note that the generating function for $q(n)$ is
\[ Q(x) = \prod_{n=1}^\infty (1+x^n) = \frac{P(x)}{P(x^2)}. \]
Hence by \eqref{Gauss}, 
\[ Q(x)  = P(x^2) \sum_{m \ge 0} x^{\frac{m(m+1)}{2}}. \]
Comparing coefficients of $x^n$ gives the required result.

By the same techniques used in \S 4.3, one finds that
\[ q(n) \sim \frac{ \e^{\frac{c}{\sqrt{2}} \sqrt{n}}}{2^23^{1/4}n^{3/4}}
\quad\text{as $n\rightarrow \infty$}. \]
It follows that the proportion of partitions with distinct parts satisfies
\begin{equation} 
\frac{q(n)}{p(n)} \sim (3n)^{1/4} \e^{-c\sqrt{n}
\left(1-\frac{1}{\sqrt{2}}\right)}
\quad\text{as $n \rightarrow \infty$}. \end{equation}

\section{Acknowledgements}
\thanks{I should like to thank
John Britnell and Atilla Mar{\'o}ti for several helpful conversations 
on the subject of partitions, and Gordon James and Peter~Neumann
for asking two of the questions which motivated this paper.
}

\def\cprime{$'$} \def\Dbar{\leavevmode\lower.6ex\hbox to 0pt{\hskip-.23ex
  \accent"16\hss}D}


\begin{thebibliography}{1}

\bibitem{HardyRamanujan}
G.~H. Hardy and S.~Ramanujan,
\newblock ``Asymptotic formulae in combinatorial analysis'',
\newblock {\em Proc. London Math. Soc. (2)}, vol. 17, pp. 75--113, 1917.

\bibitem{ErdosPartitions}
P.~Erd{\H{o}}s,
\newblock ``On an elementary proof of some asymptotic formulas in the theory of
  partitions'',
\newblock {\em Ann. of Math. (2)}, vol. 43, pp. 437--450, 1942.

\bibitem{Gauss}
C.~F. Gauss,
\newblock ``Zur {T}heorie der {T}ranscendenten {F}unctionen {G}eh{\"o}rig'',
\newblock in {\em Werke}, vol. III, pp. 435--446. K{\"o}nigslichen Gesellschaft
  der Wissenschaften zu G{\"o}ttingen, 1886.

\bibitem{Newman}
D.~J. Newman,
\newblock ``The evaluation of the constant in the formula for the number of
  partitions of {$n$}'',
\newblock {\em Amer. J. Math.}, vol. 73, pp. 599--601, 1951.

\bibitem{JK}
G.~James and A.~Kerber,
\newblock {\em The representation theory of the symmetric group}, vol.~16 of
  {\em Encyclopedia of Mathematics and its Applications},
\newblock Addison-Wesley Publishing Co., Reading, Mass., 1981.

\bibitem{Andrews}
G.~E. Andrews,
\newblock ``Two theorems of {G}auss and allied identities proved
  arithmetically'',
\newblock {\em Pacific J. Math.}, vol. 41, pp. 563--578, 1972.

\bibitem{Pak}
I.~Pak,
\newblock ``Partition bijections, a survey'',
\newblock {\em To appear in Ramanujan J.}, 69 pp.

\bibitem{MarotiLower}
A.~Mar{\'o}ti,
\newblock ``On elementary lower bounds for the partition function'',
\newblock {\em Integers}, vol. 3, pp. A10, 9 pp. (electronic), 2003.

\bibitem{KOR}
B.~K{\"u}lshammer, J.~B. Olsson, and G.~R. Robinson,
\newblock ``Generalized blocks for symmetric groups'',
\newblock {\em Invent. Math.}, vol. 151, no. 3, pp. 513--552, 2003.

\end{thebibliography}
\end{document}